\documentclass[12pt]{article}
\usepackage{amssymb}
\parskip=\medskipamount
\newtheorem{definition}{Definition}[section]
\newtheorem{lemma}[definition]{Lemma}
\newtheorem{theorem}[definition]{Theorem}
\newtheorem{proposition}[definition]{Proposition}
\newtheorem{corollary}[definition]{Corollary}
\newtheorem{remark}[definition]{Remark}

 1   
\font\ddpp=msbm10  scaled \magstep 1  

\def\QED{\hskip0.1em\hfill\null\ \null\nobreak\hfill
\kern3pt\lower1.8pt\vbox{\hrule\hbox   {\vrule\kern1pt\vbox{\kern1.7pt
\hbox{$\scriptstyle   QED$}\kern0.2pt}\kern1pt\vrule}\hrule}}
\def\R{\hbox{\ddpp R}}

\begin{document}
\title{GENERALISED CONNECTIONS OVER A VECTOR BUNDLE MAP}     
\author{Frans Cantrijn and Bavo Langerock\\
Department of Mathematical Physics and Astronomy, 
Ghent University\\
Krijgslaan 281, B-9000 Ghent, Belgium\\
e-mail: Frans.Cantrijn@rug.ac.be, Bavo.Langerock@rug.ac.be\\}
\date{\today}
\maketitle     

Short title: {\em Connections over a bundle map}
2000 MS Classification: 53C05, 58A30\\

\bigskip

\begin{abstract}
A generalised notion of connection on a fibre bundle $E$ over a manifold 
$M$ is presented. These connections are characterised by a smooth distribution 
on $E$ which projects onto a (not necessarily integrable) distribution on 
$M$ and which, in addition, is `parametrised' in some specific way by a vector 
bundle map from a prescribed vector bundle over $M$ into $TM$. Some basic properties 
of these generalised connections are investigated. Special attention is paid to
the class of linear connections over a vector bundle map. It is pointed out 
that not only the more familiar types of connections encountered in the 
literature, but also the recently studied Lie algebroid connections, can be 
recovered as special cases within this more general framework.   
\end{abstract}

\newpage
\section{Introduction}
The theory of connections undoubtedly constitutes one of the most beautiful and
most important chapters of differential geometry, which has been widely explored
in the literature (see e.g.\ \cite{GrHa,KoNo,Kol,ManSa,Sharpe,Yano}, and references
therein). Besides its purely mathematical interest, connection theory has also 
become an indispensable tool in various branches of theoretical and mathematical 
physics, as well as in other scientific disciplines which admit a proper geometric 
formulation such as, for instance, control theory and even mathematical biology
(for the latter, see \cite{Anton} for some potential applications in 
the framework of Finsler geometry).
 
Consider an arbitrary fibre bundle $\pi: E \rightarrow M$, with total space 
$E$ and base space $M$, and let $VE$ denote the canonical vertical 
distribution, i.e.\ the subbundle of $TE$ consisting of all vectors tangent 
to the fibres of $\pi$. A connection on $\pi$ (or $E$) is then given by a 
smooth distribution $HE$ on $E$, called a horizontal distribution, which is 
complementary to $VE$ and projects onto $TM$. This leads to a direct sum 
decomposition of $TE$, i.e.\ $TE = HE \oplus VE$. Note that there exist other, 
equivalent ways of characterising a connection. For instance, a connection 
on $\pi$ is sometimes defined as as a global section of the first jet bundle 
$J^1\pi$ over $E$, or also as a splitting of the short exact sequence 
\[
0 \longrightarrow VE \stackrel{i}{\longrightarrow} TE 
\stackrel{\tilde{\pi}}{\longrightarrow} 
{\pi}^{\ast}TM \longrightarrow 0, 
\]
i.e.\ a smooth map $h: {\pi}^{\ast}TM \rightarrow TE$ such that 
$\tilde{\pi} \circ h$ is the identity map on the pull-back bundle 
${\pi}^{\ast}TM$, where $i$ denotes the natural injection and 
$\tilde{\pi}$ the projection of $TE$ onto 
${\pi}^{\ast}TM$(cf.\ \cite{GrHa,ManSa,David}). 

From the above notion of connection, sometimes also called Ehresmann connection, 
one can easily derive more specific types of connections by imposing additional 
conditions on $E$ and/or $HE$. For instance, if $E$ is a 
vector bundle, it makes sense to distinguish between linear and nonlinear 
connections, depending on whether or not $HE$ is invariant under the flow of 
the canonical dilation vector field on $E$. Linear connections are often 
introduced in terms of its associated covariant derivative operator. 
If $E = TM$, it is customary to talk about a (linear or nonlinear) connection 
on $M$, instead of $TM$. In case $E$ is a principal bundle, with structure 
group $G$, and if the horizontal distribution is assumed to be $G$-invariant, 
one recovers the important notion of a principal connection. 

In the literature one can find several generalisations of the concept of  
(Ehresmann) connection introduced above, obtained by relaxing the conditions 
on $HE$. First of all, we are thinking here of the so-called 
{\it partial connections}, where the horizontal distribution $HE$ does not 
determine a full complement of $VE$. More precisely, $HE$ has zero intersection 
with $VE$, but projects onto a subbundle of $TM$, rather than onto the full 
tangent bundle (see e.g.\ \cite{KamTon}). Of special interest are 
partial connections projecting onto an integrable subbundle of $TM$, which 
play an important role in the study of the geometry of regular foliations 
(see also \cite{Kub}). 

Secondly, there also exists a notion of  {\it pseudo-connection}, introduced 
under the name of quasi-connection in a paper by Y.C. Wong \cite{Wong}. A 
fundamental role in the definition of a linear pseudo-connection on a manifold 
$M$ is played by a type $(1,1)$-tensor field on $M$ which simply becomes the 
unit tensor field in case of an ordinary linear connection. Linear 
pseudo-connections, and generalisations of it, have been studied by many 
authors (see \cite{Etayo} for a coordinate free definition of a pseudo-connection 
on a fibre bundle, and for more references to the subject). 

The inspiration for the present paper mainly stems from some recent work by 
R.L. Fernandes on a notion of `contravariant connection' in the framework
of Poisson geometry (cf.\ \cite{Fern1}). Given a Poisson manifold $(M,\Lambda)$, 
with Poisson tensor $\Lambda$ which does not have to be of constant rank, and
a principal $G$-bundle $\pi: P \rightarrow M$, a contravariant connection on $\pi$ 
is defined as a $G$-invariant bundle map $h: {\pi}^{\ast}(T^{\ast}M) 
\rightarrow TP$ over the natural vector bundle morphism 
$\sharp_{\Lambda}: T^{\ast}M \rightarrow TM$ induced by the Poisson tensor. 
This concept of connection significantly deviates from the standard one, in 
that the `horizontal' distribution $\hbox{Im}\,(h)$ may have nonzero intersection 
with the vertical subbundle $VP$ and, as for partial connections, 
projects onto a subbundle of $TM$, namely $\sharp_{\Lambda}(T^{\ast}M)$. It is 
demonstrated in \cite{Fern1} that this definition of connection leads to 
familiar concepts such as parallellism, holonomy, curvature, etc..., and, 
therefore, plays an important role in the study of global aspects of Poisson 
manifolds.  In a subsequent paper \cite{Fern2}, Fernandes has extended 
this theory by replacing the cotangent bundle of a Poisson manifold by a Lie 
algebroid over an arbitrary manifold, and the $\sharp_{\Lambda}$-map of the Poisson 
tensor by the anchor map of the Lie algebroid structure. This resulted into 
a notion of Lie algebroid connection which, in particular, 
turns out to be appropriate for studying the geometry of singular foliations. 
Fernandes' construction also covers the one given by Mackenzie \cite{Mackenzie} 
for the case of a so-called transitive Lie algebroid, where the anchor map is surjective.
   
In the present paper we will propose a general notion of connection on a fibre 
bundle $E \rightarrow M$, defined over a linear bundle morphism from an arbitrary vector 
bundle $N$ over $M$ (not necessarily a Lie algebroid) into $TM$. The relevance of the 
proposed model, in our opinion, is twofold. First of all, as will be easily recognised, 
it covers all types of connections mentioned above and, hence, it may be interesting to 
revisit some aspects of known connection theories from this broader perspective. 
Secondly, and perhaps more importantly, it may bring within the reach of connection 
theory certain geometric structures which have not yet been considered from
such a point of view.  
  
The structure of the paper is as follows. In the next section we introduce the
main definitions and describe the general framework for connections
over a vector bundle map. Section 3 is devoted to some general properties of these
connections. In Section 4 we consider various settings where this type of connections
may show up. In particular, we show how the standard Ehresmann connections, as
well as the notions of pseudo-connection, partial connection and Lie algebroid 
connection, fit into the general scheme presented here. In Section 5, special attention
is paid to the case of generalised linear connections with, among others, a discussion of
the notion of parallel transport and the construction of a suitable derivative operator.
Section 6 deals with the concepts of curvature and torsion. Generalised principal 
connections over a vector bundle map are treated rather briefly in Section 7 since, unlike
for the case of linear connections, our approach here will be very similar to the one 
adopted by Fernandes in the Lie algebroid case \cite{Fern2}. We conclude in Section 8
with some final remarks.     

{\it Notations and conventions}. The whole treatment is confined to the category of 
real, smooth (in the $C^{\infty}$ sense) geometric structures. Given a fibre bundle 
$\lambda: F \rightarrow M$, the set of all smooth sections defined on an open 
neigbourhood of a point $m \in M$ will be denoted by $\Gamma_m(\lambda)$, and we
further put $\Gamma(\lambda) = \cup_{m \in M}\Gamma_m(\lambda)$. Note, in particular,
that any global section of $\lambda$, if it exists, belongs to $\Gamma_m(\lambda)$
for all $m$. The fibre of $\lambda$ over a point $m \in M$ will be indicated by $F_m$. 
The space of smooth vector fields on a manifold $M$ will be denoted by ${\frak X}(M)$. 
Given a smooth map $f: N_1 \rightarrow N_2$ between two manifolds, we will denote the 
tangent map of $f$ by $f_{\ast}: TN_1 \rightarrow TN_2$.

\section{The general setting}
Let $M$ be a smooth $n$-dimensional manifold and $\nu: N \rightarrow M$ a vector
bundle over $M$, with $k$-dimensional fibres. Local coordinates on $M$ will be 
denoted by $(x^i)$ and the corresponding bundle coordinates on $N$ by 
$(x^i,u^{\alpha})$, with $i=1,\ldots,n$ and $\alpha=1,\ldots,k$. 
Assume we are given a vector bundle morphism $\rho: N \rightarrow TM$ over the 
identity, such that we have the following commutative diagram
 
\begin{picture}(15,4.6)(-3,4.5)
\thicklines
\put(2.3,8){$N$}
\put(6.5,8){$TM$}
\put(4.4,5.4){$M$}
\put(3,8.2){\vector(1,0){3.2}}
\put(2.7,7.6){\vector(1,-1){1.5}}
\put(6.6,7.6){\vector(-1,-1){1.5}}
\put(2.9,6.8){$\nu$}
\put(6.2,6.8){$\tau_M$}            
\put(4.4,8.4){$\rho$}
\end{picture}

with $\tau_M: TM \rightarrow M$ the canonical tangent bundle projection.
For any (local) section $s: M \rightarrow N$ of 
$\nu$, $\rho \circ s$ defines a (local) vector field on $M$. In coordinates,
$\rho$ takes the form
\begin{equation}\label{rho}
\rho(x^i,u^{\alpha}) = (x^i,{\gamma}^i_{\alpha}(x)u^{\alpha}),
\end{equation}
for some smooth functions ${\gamma}^i_{\alpha}$. For notational convenience, we
put $\hbox{Im}(\rho) = {\cal D}$. Since we do not require $\rho$ to be of constant 
rank, $\cal D$ in general will not be a vector subbundle of $TM$. Instead, it 
follows by construction that $\cal D$ determines a generalised differentiable
distribution on $M$, i.e.\ a smooth distribution in the sense of Sussmann \cite{Suss}
(see also \cite{Vais1}).

Next, let $\pi: E \rightarrow M$ be a fibre bundle over $M$, with $\ell$-dimensional
fibres and with local bundle coordinates denoted by $(x^i,y^A)$, where $i=1,\ldots,n$
and $A=1,\ldots,\ell$. We can then consider the pull-back bundle 
${\pi}^{\ast}N = \{(e,n) \in E \times N\,|\,\pi(e)=\nu(n)\}$ which
can be regarded as being fibred over $E$ as well as over $N$, with natural projections
given in coordinates by, respectively, 
\[
\tilde{\pi}_1: {\pi}^{\ast}N \longrightarrow E,\; (x^i,y^A,u^{\alpha}) \longmapsto 
(x^i,y^A)
\]
and
\[
\tilde{\pi}_2: {\pi}^{\ast}N \longrightarrow N,\; (x^i,y^A,u^{\alpha}) \longmapsto (x^i,u^{\alpha}).
\]
Note that $\tilde{\pi}_1$ is a vector bundle over $E$. In particular,
for each point $e \in E$, the fibre $(\tilde{\pi}_1)^{-1}(e)$ can be identified 
with the vector space $N_{\pi(e)}={\nu}^{-1}({\pi}(e))$. We now have all  
ingredients at hand to introduce the main concept of the present paper.

\begin{definition}\label{rhocon} {\rm A generalised connection on $\pi$ defined over 
the vector bundle morphism $\rho$, henceforth briefly called} a 
$\rho$-connection on $\pi$,{\rm is a smooth linear bundle map 
$h: {\pi}^{\ast}N \rightarrow TE$ from $\tilde{\pi}_1$ to $\tau_E$ over the identity 
on $E$, i.e.}

\begin{picture}(15,4.6)(-3,4.5)
\thicklines
\put(2.1,8){${\pi}^{\ast}N$}
\put(6.5,8){$TE$}
\put(4.4,5.4){$E$}
\put(3,8.2){\vector(1,0){3.2}}
\put(2.7,7.6){\vector(1,-1){1.5}}
\put(6.6,7.6){\vector(-1,-1){1.5}}
\put(2.9,6.8){$\tilde{\pi}_1$}
\put(6.2,6.8){$\tau_E$}            
\put(4.4,8.4){$h$}
\end{picture}

such that, in addition, the following diagram commutes:

\begin{picture}(15,4.6)(-3,4.5)
\thicklines
\put(2.4,5.4){$N$}
\put(3.2,5.6){\vector(1,0){3}}
\put(6.5,5.4){$TM$}
\put(2.7,7.6){\vector(0,-1){1.7}}
\put(6.5,8){$TE$}
\put(3.2,8.2){\vector(1,0){3}}
\put(2.2,8){${\pi}^{\ast}N$}
\put(6.8,7.6){\vector(0,-1){1.7}}
\put(4.7,5.2){$\rho$}
\put(4.7,8.4){$h$}
\put(7,6.7){${\pi}_{\ast}$}
\put(2.1,6.9){$\tilde{\pi}_2$}
\end{picture}      
\end{definition}

For a any point $(e,n) \in {\pi}^{\ast}N$, we will call $h(e,n) \in T_eE$ the 
{\it $h$-lift} of $n$ to $e$. Given any (local) section $s$ of $\nu$, we can 
define a mapping $s^h: E \rightarrow TE$ by 
\begin{equation}\label{lift}
s^h(e) = h(e, s(\pi(e))). 
\end{equation}
It is seen that, by construction, $s^h$ is smooth and verifies 
$\tau_E(s^h(e)) = e$, i.e.\ $s^h$ is a (local) vector field on $E$, 
called the {\it $h$-lift of the section $s$}. The following properties
are easily verified using the above definitions, and so we omit the proofs.

\begin{proposition}\label{sh} Given a $\rho$-connection $h$ on $\pi$, we have for 
any $s_1,s_2 \in \Gamma(\nu)$ and $f \in C^{\infty}(M)$, that:
\begin{enumerate}
\item $(s_1 + s_2)^h = s_1^h + s_2^h$;
\item $(fs)^h = (\pi^{\ast}f)s^h$;
\item $\pi_{\ast} \circ s^h = (\rho \circ s) \circ \pi$, i.e.\ the vector
fields $s^h \in {\frak X}(E)$ and $\rho \circ s \in {\frak X}(M)$ are $\pi$-related.
\end{enumerate}
\end{proposition}
   
From the definition of a $\rho$-connection $h$ it follows that for each point 
$e\in E$, the restriction of $h$ to the fibre $(\tilde{\pi}_1)^{-1}(e)$ 
of the vector bundle $\tilde{\pi}_1$, is a linear map
\[
h_e: (\tilde{\pi}_1)^{-1}(e) \cong N_{\pi(e)} \rightarrow T_eE,\; n \longmapsto h(e,n).
\]
In terms of the bundle coordinates introduced above, and taking into account 
the local expression (\ref{rho}) for $\rho$, we can write $h$ as 
\begin{equation}\label{h}
h(x^i,y^A,u^{\alpha}) = 
(x^i,y^A,{\gamma}^i_{\alpha}(x)u^{\alpha},\Gamma^A_{\alpha}(x,y)u^{\alpha}).
\end{equation}
The functions $\Gamma^A_{\alpha}$ play the role of ``connection coefficients" of
the $\rho$-connection $h$. In order to see how these functions behave under natural
coordinate transformations, take any point $(e,n) \in \pi^{\ast}N$, 
with $\pi(e) = \nu(n) = m$, and consider a change of coordinates $(x^i,y^A,u^{\alpha})
\rightarrow (\bar{x}^i,\bar{y}^A,\bar{u}^{\alpha})$ in
a neighbourhood of $(e,n)$, compatible with the underlying bundle structures:
\[
\bar{x}^i = \bar{x}^i(x), \quad  \bar{y}^A = \bar{y}^A(x,y), \quad  
\bar{u}^{\alpha} =\Lambda^{\alpha}_{\beta}(x)u^{\beta},
\]
where $\Lambda(x) = (\Lambda^{\alpha}_{\beta}(x))$ is a regular matrix. Note, first of 
all, that with respect to the bundle coordinates $(\bar{x}^i,\bar{u}^{\alpha})$ 
on $N$, the map $\rho$ can be written as $(\bar{x}^i,\bar{u}^{\alpha}) \rightarrow 
(\bar{x}^i, \bar{\gamma}^i_{\alpha}(\bar{x})\bar{u}^{\alpha})$,
with 
\[
\bar{\gamma}^i_{\alpha}(\bar{x}(x)) = \frac{\partial \bar{x}^i}{\partial x^j}(x) \gamma^j_{\beta}(x)(\Lambda^{-1})^{\beta}_{\alpha}(x)\;.
\]
Next, representing $h(e,n)$ in both coordinate systems by $(x^i,y^A,\gamma^i_{\beta}(x)u^{\beta}, \Gamma^A_{\alpha}(x,y)u^{\alpha})$ and $(\bar{x}^i,\bar{y}^A,\bar{\gamma}^i_{\beta}(\bar{x})\bar{u}^{\beta},
\bar{\Gamma}^A_{\alpha}(\bar{x},\bar{y})\bar{u}^{\alpha})$, respectively, and taking 
into account the natural coordinate transformation on $TE$, induced by the transformation 
$(x^i,y^A) \rightarrow (\bar{x}^i,\bar{y}^A)$ on $E$, one finds after a tedious, but 
straightforward computation, the following transformation law for the connection
coefficients associated to a general $\rho$-connection:
\begin{equation}\label{Gamma}
\bar{\Gamma}^A_{\alpha}(\bar{x}(x),\bar{y}(x,y)) = 
\left( \frac{\partial \bar{y}^A}{\partial x^j}(x,y)\gamma^j_{\beta}(x)
+ \frac{\partial \bar{y}^A}{\partial y^B}(x,y)\Gamma^B_{\beta}(x,y)\right)
(\Lambda^{-1})^{\beta}_{\alpha}(x)\;.
\end{equation}

Henceforth, given a $\rho$-connection $h$ we will put for brevity: 
$\hbox{Im}$$(h) = \cal Q$. This determines a smooth generalised distribution on 
$E$ which projects onto ${\cal D}\;(= \hbox{Im}(\rho))$. We refrain from calling 
$\cal Q$ a horizontal distribution since for arbitrary $e \in E$ it may be 
that ${\cal Q}_e$ has non-zero intersection with $V_eE$. Moreover, in general 
${\cal Q}_e + V_eE \neq T_eE$, i.e.\ ${\cal Q}_e$ and $V_eE$ do not necessarily span 
the full tangent space $T_eE$.

In the above, $\pi:E \rightarrow M$ always represented an arbitrary fibre bundle over 
$M$. Some interesting types of $\rho$-connections are obtained when imposing 
additional conditions on $E$. First of all, if $\pi: E=P \rightarrow M$ is a 
principal $G$-bundle with a, say, right Lie group action 
$\Phi: P \times G \rightarrow P,\; (e,g) \mapsto (\Phi(e,g)=)\Phi_g(e)=eg$, then 
a $\rho$-connection $h$ on $\pi$ is called a {\it principal $\rho$-connection\/} if
\[
(\Phi_g)_{\ast}(h(e,n))=h(eg,n)
\]
for all $g \in G$ and $(e,n) \in {\pi}^{\ast}N$. In particular, this implies that 
the associated distribution $\cal Q$ is $G$-invariant.

Next, assume $E$ is the total space of a vector bundle over $M$. Then, 
$\tilde{\pi}_2: {\pi}^{\ast}N \rightarrow N$ is also a vector bundle, the
fibres of which can be identified with those of $\pi$. Let $\{\phi_t\}$ represent
the flow of the canonical dilation vector field on $E$ , i.e.\ in natural vector 
bundle coordinates $(x^i,y^A)$ on $E$ we have $\phi_t(x^i,y^A) = (x^i, e^ty^A)$.
We will then say that $h$ is a {\it linear $\rho$-connection\/} if for all $(e,n) \in
\pi^{\ast}N$
\[
(\phi_t)_{\ast}(h(e,n)) = h(\phi_t(e),n)\,.
\]
It is not difficult to see that this implies that the connection coefficients,
appearing in (\ref{h}), are of the form $\Gamma^A_{\alpha}(x,y) = 
\Gamma^A_{{\alpha}B}(x)y^B$.

Let us return to the general situation described by Definition \ref{rhocon},
with $\pi: E \rightarrow M$ an arbitrary fibre bundle over $M$. 
Regarding $TE$ as a vector bundle over $TM$, with projection $\pi_{\ast}$, 
we can define the pull-back bundle 
$\rho^{\ast}TE = \{(n,w) \in N \times TE\;|\; \rho(n) = \pi_{\ast}(w)\}$. Clearly,
if $(n,w) \in \rho^{\ast}TE$, with $\tau_E(w) = e$, then $(e,n) \in \pi^{\ast}N$
and, given a $\rho$-connection $h$ on $\pi$, one easily verifies that
\[
\pi_{\ast}(w - h(e,n)) = 0.
\]
Hence, one can define a mapping $V: \rho^{\ast}TE \rightarrow VE$ by
\begin{equation}\label{V}
V(n,w) = w - h(e,n) \quad \hbox{with}\; e = \tau_E(w).
\end{equation}
Note that $\rho^{\ast}TE$ admits the structure of a vector bundle over $E$, with 
fibre over $e \in E$ given by $N_m \times T_eE$, where $m = \pi(e)$. With respect 
to this stucture, it is straightforward to check that $V$ is a vector bundle 
morphism over the identity on $E$. If, in the appropriate bundle coordinates, 
$(x^i,u^{\alpha})$ are the coordinates of a point $n \in N$ and 
$(\bar{x}^i,y^A,v^i,w^A)$ those of $w \in TE$, then the condition that $(n,w)$ 
represents an element of $\rho^{\ast}TE$ boils down to the requirement that 
$x^i = \bar{x}^i$ and $v^i = \gamma^i_{\alpha}(x)u^{\alpha}$. Therefore, natural
bundle coordinates on $\rho^{\ast}TE$, induced by those on $N$ and $TE$, are $(x^i,u^{\alpha},y^A,w^A)$. In terms of the latter, the mapping $V$ can now be written as
\[
V(x^i,u^{\alpha},y^A,w^A) = (x^i,y^A,0,w^A - \Gamma^A_{\alpha}(x,y)u^{\alpha}), 
\]
where $\Gamma^A_{\alpha}(x,y)$ are the connection coefficients associated to $h$.

In case $\pi: E \rightarrow M$ is a vector bundle, it is well-known that there 
exists a canonical isomorphism between $VE$ and the fibred product 
$E \times_ME \;(\cong \pi^{\ast}E)$. Denote by 
$p_2: VE \cong E \times_M E \rightarrow E$ the projection onto the second factor,
i.e.\ in coordinates: $p_2(x^i,y^A,0,w^A) = (x^i,w^A)$. Given a (not necessarily linear)
$\rho$-connection $h$ on $\pi$, we can define a mapping 
$K: \rho^{\ast}TE \rightarrow E$ by 
\begin{equation}\label{K1}
K (n,w) = (p_2 \circ V)(n,w) \quad \hbox{for all} \quad (n,w)\in \rho^{\ast}TE.
\end{equation}
In coordinates, taking into account the above expression for $V$, this reads
\begin{equation}\label{K2}
K(x^i,u^{\alpha},y^A,w^A) = (x^i, w^A - \Gamma^A_{\alpha}(x,y)u^{\alpha}).
\end{equation}
The mapping $K$ will be called the {\em connection map} (associated to the given
$\rho$-connection), in analogy with the connection map associated 
to an ordinary connection on a vector bundle (see e.g.\ \cite{Vilms}).

To close this section, we now introduce a special class of curves in $N$ which 
will play a central role, among others, when considering a notion of 
parallel transport in the framework of generalised connections over a vector 
bundle map. By a smooth curve in a manifold $Q$ we will always mean
a $C^{\infty}$ map $c: I \rightarrow Q$, where $I \subseteq \R$ may be either an 
open or a closed (compact) interval. In the latter case, the denominations ``path" 
or ``arc" are also frequently used in the literature but, for simplicity, we will 
make no distinction in terminology between both cases. For a curve defined on
a closed interval, say $[0,1]$, it is tacitly assumed that it admits a smooth
extension to an open interval containing $[0,1]$.     

For a given curve $c$ in $N$ we put $\tilde{c} = \nu \circ c$, i.e.\ $\tilde{c}$ 
is the projection of $c$ onto $M$.

\begin{definition}\label{curve} {\rm A smooth curve $c: I \rightarrow N$, is called 
a} $\rho$-admissible curve {\rm if 
\[
(\rho \circ c)(t) = \dot{\tilde{c}}(t)\,,
\]
for all $t \in I$. The projection $\tilde{c}$ of a $\rho$-admissible curve will be 
called a} base curve.
\end{definition}

For a detailed treatment of certain aspects of the geometry related to
connections over a bundle map, it will be necessary to extend the notion of
$\rho$-admissibility to allow for continuous, piecewise 
smooth curves, and even for curves admitting a finite number of discontinuities, 
which are such that the projections of these curves onto $M$ are piecewise smooth.
Such curves will then also be called $\rho$-admissible, provided each of
its `smooth components' is $\rho$-admissible. It should also be pointed out that, 
in principle, weaker types of smoothness (e.g. $C^1$) would 
have been sufficient for most considerations. For the purpose of the present paper,
however, we will confine ourselves to the class of smooth $C^{\infty}$-curves.

Occasionally, if no confusion can arrise, we will also simply refer to a 
$\rho$-admissible curve $c$ as an ``admissible curve". In coordinates, if we put 
$c(t) = (x^i(t),u^{\alpha}(t))$, the condition for $\rho$-admissibility reads
\[
\dot{x}^i(t) = \gamma^i_{\alpha}(x(t))u^{\alpha}(t).
\]
From the definition it immediately follows that a base curve is 
everywhere tangent to the (generalised) distribution $\cal D$. In particular, for each 
section $s$ of $\nu$, the integral curves of the vector field $\rho \circ s
(= \rho(s))$ are base curves. If $\cal D$ is an integrable distribution, it 
is follows that a smooth base curve is contained in a leaf of the induced foliation 
of $M$. 

An important observation is that, due to the fact that $\rho$ need 
not be injective, there may be different $\rho$-admissible curves passing through a 
given point $n \in N$ which project onto the same base curve.  
Note also that a curve $c$ in $N$ whose image belongs to $\ker (\rho)$, will be 
$\rho$-admissible iff $c$ is contained in a fibre of $N$, and $\tilde{c}$ then
reduces to a point. Finally, if a smooth base curve $\tilde{c}: I \rightarrow M$ is an
immersion, i.e.\ $\dot{\tilde{c}}(t) \neq 0$ for all $t \in I$, it follows that the 
image of any $\rho$-admissible curve which projects onto $\tilde{c}$ must have empty intersection with $\ker (\rho)$.

\section{Some general properties}

As observed above, the distribution ${\cal Q}$ defined by a $\rho$-connection $h$ on
a fibre bundle $\pi: E \rightarrow M$, in general may have nonzero intersection with 
the vertical subbundle $VE$ of $TE$. The extent by which $\cal{Q}$ fails to be a (full) complement of $VE$ is characterised by the following proposition.

\begin{proposition} For any $m \in M$ and $e \in E_m$ we have\\
\begin{equation}\label{intersection}
{\cal Q}_e \cap V_eE \cong \ker({\rho}_m)/\ker(h_e),
\end{equation}
(where ${\rho}_m$ and $h_e$ are the linear maps induced by the restrictions
of $\rho$ and $h$, respectively, to the fibre $N_m$ of $N$), and   
\begin{equation}\label{sum}
{\cal Q}_e + V_eE = T_eE \quad \iff \quad {\cal D}_m = T_mM.
\end{equation}
\end{proposition}
{\bf Proof.} For $w \in T_eE$, with $\pi(e) = m$, we immediately have that 
$w \in {\cal Q}_e \cap V_eE$ iff $w = h(e,n) = h_e(n)$ for some $n \in N_m$, and
$0 = {\pi}_{\ast}(w) = {\pi}_{\ast}(h(e,n)) = \rho(n) = {\rho}_m(n)$. Hence,
\[
w \in {\cal Q}_e \cap V_eE \iff w \in h_e(\ker(\rho_m)).
\]  
From the definition of $h$ one can deduce that $\ker(h_e) \subset \ker({\rho}_m)$ 
and it then readily follows that $h_e(\ker({\rho}_m)) \cong \ker({\rho}_m)/\ker(h_e)$,
which completes the proof of (\ref{intersection}).

Next, assume that ${\cal Q}_e + V_eE = T_eE$, for $e \in E_m$. For any $v \in T_mM$
one can always find a $w \in T_eE$ such that ${\pi}_{\ast}(w) = v$. The given 
assumption implies that $w$ can be written as $w = h_e(n) + \tilde{w}$, 
for some $n \in N_m$ and $\tilde{w} \in V_eE$, and this, in turn, gives
\[
v = {\pi}_{\ast}(w) = {\pi}_{\ast}(h_e(n)) = \rho_m(n),
\]
i.e.\ $v \in \hbox{Im}(\rho_m) = {\cal D}_m$. Since $v \in T_mM$ was chosen arbitrarily,
this proves that ${\cal D}_m = T_mM$. Conversely, assume ${\cal D}_m = T_mM$. For
any $w \in T_eE$ we then have that ${\pi}_{\ast}(w) = \rho_m(n)$ for some $n \in N_m$,
from which it follows that ${\pi}_{\ast}(w - h_e(n)) = {\pi}_{\ast}(w) - \rho_m(n) = 0$,
and so $w - h_e(n) \in V_eE$. This completes the proof of the equivalence (\ref{sum}).
\QED 

From this proposition one can readily deduce the following result.

\begin{corollary} The distribution $\cal Q$ defines a genuine (Ehresmann) connection 
on $\pi$ iff $\rho(N) = TM$ and $\ker(\rho_m) = \ker(h_e)$ for all $m \in M$ and 
$e \in E_m$.
\end{corollary}

Whereas a $\rho$-connection $h$ determines a (generalised) distribution $\cal Q$ on $E$
which projects onto $\cal D$, the converse is certainly not true in general. Moreover,
if a distribution $\cal Q$ can be associated to a $\rho$-connection, the latter
need not be uniquely determined. A sufficient condition for a distribution on $E$
to correspond to a unique $\rho$-connection is that it determines a (not necessarily 
full) complement of $VE$.

\begin{proposition} Let $\cal Q$ be a smooth generalised distribution on $E$ such that
(i) $\pi_*({\cal Q}) = {\cal D}$, and (ii) ${\cal Q}_e \cap V_eE = \{0\}$ for all 
$e \in E$, then there exists a unique $\rho$-connection $h$ such that
${\cal Q} =$ Im($h$).
\end{proposition} 
{\bf Proof.} For each point $e \in E$, we can construct a map 
$h_e: N_m \rightarrow T_eE$, where $m = \pi(e)$, by putting 
\[
\{h_e(n)\} = {\cal Q}_e \cap \left(({\pi}_*)_{|T_eE}\right)^{-1}(\rho_m(n)),
\]
for all $n \in N_m$. From the given assumptions (i) and (ii), it follows that 
this prescription uniquely determines a point $h_e(n)$. Furthermore, using some simple
set-theoretic arguments, it is not difficult to verify that the resulting map 
$h_e$ is linear. Next, we can `glue' these linear maps together to a smooth bundle map 
$h: \pi^*N \rightarrow TE$ with $h(e,n) = h_e(n)$. It is then straightforward to see
that, by construction, $h$ verifies all properties of a $\rho$-connection. 

Finally, uniqueness of $h$ can be proved as follows. Let $h': \pi^{\ast}N 
\rightarrow TE$ be another $\rho$-connection for which $\hbox{Im}(h')= {\cal Q}$. 
Then, for each $(e,n) \in \pi^{\ast}N$, with $\pi(e) = \nu(n) = m$, there exists 
a $n' \in N_m$ such that $h(e,n) = h'(e,n')$. The definition of a $\rho$-connection 
then implies that $\rho(n) = \rho(n')$. Now, from (3) and the assumption (ii) it follows
that $\ker(\rho_m)=\ker(h'_e)$ and, hence, $h'(e,n) = h'(e,n') = h(e,n)$, which indeed
proves uniqueness of the $\rho$-connection.    
\QED

Herewith we can now prove the following result.

\begin{theorem} Given a vector bundle $\nu: N \rightarrow M$, a vector bundle
morphism $\rho: N \rightarrow TM$ such that $\nu = \tau_M \circ \rho$, and a fibre 
bundle $\pi: E \rightarrow M$. Then, there always exists a $\rho$-connection on 
$\pi$.
\end{theorem}
{\bf Proof.} The proof immediately follows from the previous proposition and the 
well-known property that on each fibre bundle one can always construct an ordinary 
connection (see e.g.\ \cite{ManSa}). Indeed, take an arbitrary connection on
$\pi$ with horizontal distribution denoted by $HE$, such that $TE = HE \oplus VE$.
Then, putting ${\cal Q} = (\pi_{\ast})^{-1}({\cal D}) \cap HE$, it is easily verified 
that $\cal Q$ defines a (generalised) distribution on $E$, satisfying the conditions
of Proposition 3.3.
\QED

Note that the $\rho$-connections referred to in Proposition 3.3 and, consequently, 
also the one constructed in the previous theorem, are of a special type in the sense 
that the corresponding distribution $\cal Q$ is `transverse' to $VE$,
i.e.\ ${\cal Q}_e \cap V_eE = \{0\}$ for all $e \in E$. With a slight abuse of 
terminology, we will call such a $\rho$-connection a {\it partial connection} on $\pi$. 
If the distribution $\cal Q$ has constant rank it determines indeed a partial connection
in the ordinary sense (see the Introduction).

\begin{remark} {\rm The notion of partial connection, as defined above, also 
corresponds to (and reduces to) what Fernandes has called ${\cal F}$-connections 
in his treatment of contravariant connections on Poisson manifolds and connections 
on Lie algebroids \cite{Fern1,Fern2}.}
\end{remark}

Assume now that $\rho$ has constant rank. Then, $\hbox{Im}(\rho)$ is a 
vector subbundle of $TM$, with canonical injection 
$i:\hbox{Im}(\rho) \hookrightarrow TM$. 

\begin{proposition} If $\rho$ has constant rank, then for every $\rho$-connection 
$h$ on a fibre bundle $\pi: E \rightarrow M$ there is a $i$-connection 
$\bar{h}$ on $\pi$ such that $\hbox{Im}$$(h) = \hbox{Im}$$(\bar{h})$ 
iff $h$ is a partial connection.
\end{proposition} 
{\bf Proof.} If $h$ is a partial connection, we know from the above that
$\ker(\rho_m) = \ker(h_e)$ for all $m \in M$ and $e \in E_m$.
We can then define a mapping $\bar{h}: \pi^{\ast}\hbox{Im}(\rho) \rightarrow TE$
as follows: for $n \in N_m$ and $e \in E_m$, put
\[
\bar{h}(e,\rho(n)) = h(e,n)\;.
\]
From the fact that $h$ is a partial connection it follows that $\bar{h}$ is
well defined, and it is straightforward to check that it is a generalised
connection over $i$, determining the same distribution on $E$ as $h$.

Conversely, assume that there exists a $i$-connection $\bar{h}$ on $\pi$, having
the same image as a given $\rho$-connection $h$. In particular, this implies
that for all $(e,n) \in \pi^{\ast}N$, with $\nu(n) = \pi(e) = m$, there exists 
a $n' \in N_m$ such that $h(e,n) = \bar{h}(e,\rho(n'))$. Since, obviously,
$\ker (i_{\rho(n')}) = 0$, we also have $\ker (\bar{h}_e) = 0$, from which one 
can readily deduce that $\ker (h_e) = \ker (\rho_m)$ and, hence, $h$ is a partial 
connection.
\QED

Next, consider the case where $\cal D (= \hbox{Im}(\rho))$ is a (generalised) 
integrable distribution, inducing a foliation of $M$, i.e.: through 
each point of $M$ passes a maximal integral manifold of $\cal D$, called
a leaf of the foliation. These leaves are immersed submanifolds of $M$ which 
need not all have the same dimension since $\rho$ (and, therefore, also $\cal D$) 
is not assumed here to be of constant rank. 
Let $S$ be an arbitrary leaf of the foliation and let $i_S: S \hookrightarrow M$ 
denote the natural injection. In particular, $i_S$ is an injective immersion.
The pull-back bundle $i_S^{\ast}N$ of $\nu: N \rightarrow M$ by $i_S$ is a vector 
bundle over $S$ (which can be identified with the restriction $N_{|S}$). 
Since ${i_S}_{\ast}$ is injective and since for each
$m \in S$, $T_mS = \rho(N_m)$, we can define a vector bundle morphism 
$\rho_S:i_S^{\ast}N \rightarrow TS$ in an implicit way by the following
prescription: for each $(m,n) \in i^{\ast}_SN$,
\[
{i_S}_{\ast}(\rho_S(m,n)) = \rho(n) \;.
\]
One can then show that a $\rho$-connection always induces a $\rho_S$-connection.
(For the analogous result in the case of connections on Lie algebroids, 
see \cite{Fern2}.) 

\begin{proposition} Let $\hbox{Im}(\rho)$ be an integrable distribution and
$S$ a leaf of the corresponding foliation of $M$. Then, every $\rho$-connection on 
a fibre bundle $\pi: E \rightarrow M$ induces a $\rho_S$-connection on 
the pull-back bundle $\pi_S: i_S^{\ast}E \rightarrow S$.
\end{proposition}
{\bf Proof.} Let $h$ be a $\rho$-connection on $\pi$. Consider the pull-back
bundle $\pi_S^{\ast}(i_S^{\ast}N)$, admitting the double fibration
$(\widetilde{\pi_S})_1: \pi_S^{\ast}(i_S^{\ast}N) \rightarrow i_S^{\ast}E$
and $(\widetilde{\pi_S})_2: \pi_S^{\ast}(i_S^{\ast}N) \rightarrow i_S^{\ast}N$.
An element of $\pi_S^{\ast}(i_S^{\ast}N)$ can be identified with a triple $(m,e,n)$, 
with $m \in S, e \in E_m, n \in N_m$. Now, define the mapping
$h_S: \pi_S^{\ast}(i_S^{\ast}N) \longrightarrow Ti_S^{\ast}E$
by
\begin{equation}\label{hS}
h_S(m,e,n):= (\rho_S(m,n),h(e,n))\;,
\end{equation}
where on the right-hand side we have used the canonical identification 
$Ti_S^{\ast}E \cong ({i_S}_{\ast})^{\ast}TE$. Note that, as such, $h_S$ is 
well defined since ${i_S}_{\ast}(\rho_S(m,n)) = \rho(n) = \pi_{\ast}(h(e,n))$.
It is then easily verified that $h_S$ is a linear bundle morphism from
$(\widetilde{\pi_S})_1$ to $\tau_{i_S^{\ast}E}$, satisfying
${\pi_S}_{\ast}(h_S(m,e,n)) = \rho_S(m,n)$.
\QED   

Under the assumptions of the previous proposition, let us put 
$\hbox{Im}$$(h_S) = {\cal Q}_S$ and, as before, $\hbox{Im}$$(h) = {\cal Q}$.
Put 
\[
\varphi_S: i_S^{\ast}E \rightarrow E,\; (m,e) \mapsto e\;,
\] 
such that $\pi \circ \varphi_S = i_S \circ \pi_S$. This is an immersion and 
we clearly have that ${\varphi_S}_{\ast}({\cal Q}_S)_{(m,e)} = {\cal Q}_e$
for all $(m,e) \in i_S^{\ast}E$. The following corollary shows that 
in case $h$ is a partial connection, this property uniquely characterises the 
$\rho_S$-connection $h_S$. 
 
\begin{corollary} If $h$ is a partial connection, then $h_S$, defined by (\ref{hS}), 
is the unique $\rho_S$-connection satisfying 
${\varphi_S}_{\ast}(\hbox{Im}$$(h_S))_{(m,e)} = {\cal Q}_e$ for
all $(m,e) \in i_S^{\ast}E$.
\end{corollary}
{\bf Proof.} First, recall that $h$ being a partial connection means that 
$\ker (\rho_m) = \ker (h_e)$ for all $m \in M, e \in E_m$. Let $\hat{h}: \pi_S^{\ast}(i_S^{\ast}N) \longrightarrow Ti_S^{\ast}E$ be any 
$\rho_S$-connection such that ${\varphi_S}_{\ast}(\hbox{Im}(\hat{h}))_{(m,e)} 
= {\cal Q}_e$ for all $(m,e) \in i_S^{\ast}E$. This implies
that for any $(m,e,n) \in \pi_S^{\ast}(i_S^{\ast}N)$ there exists
a $n' \in N_m$ such that ${\varphi_S}_{\ast}(\hat{h}(m,e,n))= h(e,n')$.  
On the other hand, from (\ref{hS}) we derive that 
${\varphi_S}_{\ast}(h_S(m,e,n'))= h(e,n')$. Hence, 
$h_S(m,e,n') - \hat{h}(m,e,n) \in \ker {\varphi_S}_{\ast}$. But $\varphi_S$
is an immersion, hence $h_S(m,e,n') = \hat{h}(m,e,n)$. From the definition
of a $\rho_S$-connection it then follows that $\rho_S(m,n) = \rho_S(m,n')$,
which implies $\rho(n) = \rho(n')$, i.e.\ $n - n' \in \ker(\rho_m)$.
From the assumption that $h$ is a partial connection we then deduce
that ${\varphi_S}_{\ast}(\hat{h}(m,e,n))= h(e,n') = h(e,n) = 
{\varphi_S}_{\ast}(h_S(m,e,n))$ which, again in view of the injectivity of 
${\varphi_S}_{\ast}$, finally shows that $\hat{h}(m,e,n) = h_S(m,e,n)$ for
all $(m,e,n)\in \pi_S^{\ast}(i_S^{\ast}N)$.
\QED

In the next section we shall describe several situations where generalised connections 
over a vector bundle map may be considered. In particular, we will see how the various 
types of connections mentioned in the Introduction can be recovered as special cases
of the general notion of connection put forward in Definition \ref{rhocon}.

\section{Special cases}
 
(i) If we put $N=TM$, $\nu = \tau_M$ and $\rho = \hbox{Id}_{TM}$ 
(the identity map on $TM$), Definition \ref{rhocon} reduces to that of an ordinary 
connection (an Ehresmann connection) on $\pi$, with 
$h: {\pi}^{\ast}TM \rightarrow TE$ defining a splitting of the short exact sequence 
$0 \rightarrow VE \rightarrow TE \rightarrow {\pi}^{\ast}TM \rightarrow 0$ and
$\hbox{Im}(h) = HE$ the horizontal distribution of the connection. In particular,
for $E = TM$ we recover the standard notion of (linear or nonlinear) connection on 
a manifold $M$(see also \cite{Vilms}).

(ii) Let $N$ be a subbundle of $TM$, $\nu = (\tau_M)_{|N}$, and 
$\rho = i_N: N \hookrightarrow TM$ the canonical injection. In this case, each 
$\rho$-connection $h$ on a fibre bundle $\pi: E \rightarrow M$ is a partial 
connection. Indeed, since for all $m \in M$ we have $\ker ((i_N)_m) = \{0\}$, it
follows from (\ref{intersection}) that ${\cal Q}_e \cap V_eE = \{0\}$ for all $e \in E$.
Moreover, $h$ is now necessarily injective, implying that $\cal Q$ is a constant rank
distribution and, therefore, we are dealing with a partial connection in the ordinary
sense. Partial connections are considered in particular in those
cases where $N$ defines an a regular integrable distribution on $M$ 
(see e.g.\ \cite{KamTon}). The horizontal subspaces ${\cal Q}_e$ then project onto
the tangent spaces to the leaves of the induced foliation. 
But partial connections also make their appearance, for instance, in the framework 
of sub-Riemannian geometry, where $N$ is a subbundle of $TM$ equipped with a 
nondegenerate bundle metric (see e.g.\ \cite{El}).
 
(iii) If $\nu: N \rightarrow M$ is a Lie algebroid over $M$, with anchor map $\rho$, 
we recover the notion of {\it Lie algebroid connection} studied by 
Fernandes \cite{Fern2}. By definition of a Lie algebroid, the anchor map induces a Lie 
algebra morphism from the Lie algebra of sections of $\nu$ into the Lie algebra of 
vector fields on $M$. Consequently, in this case $\hbox{Im}(\rho) = {\cal D}$ is an 
involutive generalised distribution, determining a (possibly singular) foliation 
$\cal F$ of $M$. Given a $\rho$-connection $h$ on a fibre bundle 
$\pi: E \rightarrow M$, with associated distribution ${\cal Q}$, we have that for 
each $e \in E$ the subspace ${\cal Q}_e$ of $T_eE$ projects onto the tangent space 
at $\pi(e)$ to the leaf of $\cal F$ passing through $\pi(e)$. Here, unlike the case 
of a partial connection, ${\cal Q}$ may have a nonzero intersection with the vertical distribution $VE$.      

A particular instance of a Lie algebroid is obtained when $M$ admits a 
Poisson structure, with Poisson tensor $\Lambda$, and $N = T^{\ast}M$. The anchor
map $\rho$ is then given by the natural vector bundle morphism induced by 
$\Lambda$, i.e.\ 
$\sharp_{\Lambda}: T^{\ast}M \rightarrow TM, \; \alpha_m \mapsto \Lambda_m(\alpha_m,.)$.
This case was also studied extensively by Fernandes \cite{Fern1}. Connections over $\sharp_{\Lambda}$ were then called {\it contravariant connections}, following 
I. Vaisman who introduced a notion of contravariant derivative in the framework of the 
geometric quantisation of Poisson manifolds \cite{Vais}.  

(iv) Let again $N=TM$, $\nu = \tau_M$ and let $\rho$ be the tangent bundle morphism 
induced by a type $(1,1)$-tensor field $A$ on $M$. A $\rho$-connection 
on then corresponds to what is also known as a {\it pseudo-connection} 
with fundamental tensor field $A$ (cf. \cite{Etayo,Wong}). 

Consider the case where $A$ has vanishing Nijenhuis torsion, i.e.\ 
${\cal N}_A = 0$, with ${\cal N}_A$ the type $(1,2)$-tensor field defined by
$1/2{\cal N}_A(X,Y) = A^2([X,Y]) + [A(X),A(Y)] - A([A(X),Y]) - A([X,A(Y)])$ for
arbitrary $X,Y \in {\frak X}(M)$. The pair $(M,A)$ is sometimes called a {\it Nijenhuis 
manifold}, with {\it Nijenhuis tensor} $A$.
One may then define a new bracket on ${\frak X}(M)$ according to
\begin{equation}\label{Nijenhuis}
[X,Y]_A := [A(X),Y] + [X,A(Y)] - A([X,Y])\,.
\end{equation}
Using the fact that ${\cal N}_A = 0$, it follows after some tedious but straightforward computations that $[\;,\;]_A$ is again a Lie bracket on ${\frak X}(M)$ and that, moreover, 
$A([X,Y]_A) = [A(X),A(Y)]$ and $[X,fY]_A = f[X,Y]_A + A(X)(f)Y$ for all 
$X,Y \in {\frak X}(M)$ and $f \in C^{\infty}(M)$ (see e.g.\ \cite{Kos}). 
Consequently, $TM$ becomes a Lie algebroid over $M$ with bracket $[\;,\;]_A$ and 
anchor map $A$ (regarded as a bundle map from $TM$ into itself), and a 
pseudo-connection whose fundamental tensor field $A$ is a Nijenhuis tensor, is a Lie 
algebroid connection.

(v) An immediate extension of the previous case is obtained when considering 
an arbitrary vector valued tensor field ${\cal K} \in {\cal T}^r_s(M) \otimes 
{\frak X}(M)$ on $M$, where ${\cal T}^r_s(M)$ denotes
the $C^{\infty}(M)$-module of smooth type $(r,s)$-tensor fields, i.e.\ tensor
fields of contravariant order $r$ and covariant order $s$. 
Putting $N = T^s_r(TM)$, the vector bundle of type $(s,r)$-tensors on $M$, and 
$\rho: T^s_r(M) \rightarrow TM$ the natural bundle 
morphism over $M$ induced by $\cal K$, i.e.
\[
\rho(v_1 \otimes \ldots \otimes v_s \otimes \alpha_1 \otimes \ldots \otimes \alpha_r)
= {\cal K}(v_1, \ldots, v_s; \alpha_1, \ldots, \alpha_r),
\]
for arbitrary $x \in M$, $v_i \in T_xM$ and $\alpha_j \in T^{\ast}_xM$,
then one can consider $\rho$-connections on a fibre bundle $E$
over $M$ as connections which, in some sense, are ``parametrised" by $(s,r)$-tensors. 
Clearly, the pseudo-connections mentioned above, as well as the contravariant (Poisson) connections, belong to this category.

(vi) Another example, which also fits into the previous category, is provided by
sub-Riemannian geometry. A sub-Riemannian structure consists of a triple $(M,Q,g)$,
where $M$ is a smooth manifold, $Q$ a distribution on $M$ of constant rank (i.e.\ a
vector subbundle of $TM$) and $g$ a positive definite bundle metric on $Q$ (see
e.g.\ \cite{El,Strich}). Herewith one can associate a vector bundle morphism
$\sharp_{g}: T^{\ast}M \rightarrow Q$ which is uniquely determined by
\[
g(v_m,\sharp_g(\alpha_m)) = \langle v_m,\alpha_m \rangle \;,
\]
for all $v_m \in Q_m$, and with $\langle \;,\;\rangle$ denoting the natural pairing
between $T_mM$ and $T^{\ast}_mM$. One can easily verify that $\ker (\sharp_g) = Q^o$, 
the annihilator of $Q$ in $T^{\ast}M$. Since $\sharp_g$ can also be regarded as a
smooth bundle morphism over the identity from $T^{\ast}M$ into $TM$, we may thus look 
for connections over the vector bundle map $\sharp_g$ in the sense of Definition 
\ref{rhocon}(with $N = T^{\ast}M$ and $\rho = \sharp_g$). Such connections will be
considered in a forthcoming paper (\cite{Bavo1}).   

\section{Linear $\rho$-connections}
In this section we assume that $\pi: E \rightarrow M$ is a vector bundle and that
$h: \pi^{\ast}N \rightarrow TE$ defines a linear $\rho$-connection on $\pi$ 
(cf.\ Section 2). Recall that, in terms of natural bundle coordinates 
$(x^i,u^{\alpha})$ and $(x^i,y^A)$ on $N$ and $E$, respectively, and with 
$\rho$ given by (\ref{rho}), the bundle map $h$ is of the form
\[
h(x^i,y^A,u^{\alpha}) = 
(x^i,y^A,\gamma^i_{\alpha}(x)u^{\alpha},\Gamma^A_{\alpha B}(x)u^{\alpha}y^B)\,.
\]
Considering an admissible coordinate transformation in a neighbourhood of some point
$(e,n) \in \pi^{\ast}N$, of the form
\[
\bar{x}^i = \bar{x}^i(x), \quad  \bar{y}^A = \Xi^A_B(x)y^B, \quad  
\bar{u}^{\alpha} =\Lambda^{\alpha}_{\beta}(x)u^{\beta},
\]
where $\Xi(x) = (\Xi^A_B(x))$ and $\Lambda(x)= (\Lambda^{\alpha}_{\beta}(x))$
are regular matrices, it can be easily deduced from the general tranformation
law (\ref{Gamma}) for the connection coefficients $\Gamma^A_{\alpha}$, that
the $\Gamma^A_{\alpha B}$ transform according to
\[
\bar{\Gamma}^A_{\alpha B}(\bar{x}(x)) = 
\left(\frac{\partial \Xi^A_C}{\partial x^k}(x)\gamma^k_{\beta}(x)
+ \Gamma^D_{\beta C}(x)\Xi^A_D(x)\right)
(\Xi^{-1})^C_B(x)(\Lambda^{-1})^{\beta}_{\alpha}(x).
\]
In particular, if $N = E = TM$ with $\nu = \pi = \tau_M$, and $\rho = \hbox{Id}_{TM}$, 
we have that both $\Lambda(x)$ and $\Xi(x)$ reduce to the Jacobian matrix 
$(\partial \bar{x}^i/\partial x^j)$ of the coordinate transformation on the base 
manifold $M$, and we recover the standard transformation law for the connection 
coefficients (``Christoffel symbols") of a linear connection on a manifold.

\subsection{Parallel transport}
We now aim at defining a notion of parallel transport for linear $\rho$-connections. 
In the next proposition, we first show that a $\rho$-admissible curve on $N$ 
(cf.\ Definition \ref{curve}) can always be lifted to a curve on $E$ which is 
everywhere tangent to the generalised distribution $\hbox{Im}(h) = \cal Q$ determined 
by the given linear $\rho$-connection.

\begin{proposition}\label{liftcurve} Consider a smooth $\rho$-admissible curve 
$c: [0,1] \rightarrow N$, with $c(0)=n_0$. Then, for each $e_0 \in E_{\nu(n_0)}$,
there exists a uniquely defined curve $c^h: [0,1] \rightarrow E$ such that 
$c^h(0)=e_0$, $(\pi \circ c^h)(t) = (\nu \circ c)(t)$ for all $t \in [0,1]$, and
\[
\dot{c}^h (t) = h(c^h(t),c(t))\,.
\]
\end{proposition} 
{\bf Proof.} The proof proceeds along the same lines as for the construction
of the horizontal lift of curves in standard connection theory. First, consider
a coordinate neighbourhood $U \subset M$ which is locally trivialising with
respect to both vector bundle structures $\nu$ and $\pi$. Coordinates on 
$\nu^{-1}(U)$ and on $\pi^{-1}(U)$ are denoted by $(x^i,u^{\alpha})$ and $(x^i,y^A)$, respectively. Assume now that the image of the given $\rho$-admissible curve $c$
is contained in $\nu^{-1}(U)$, with $c(0) = n_0 = (x^i_0,u^{\alpha}_0)$. Then, putting 
$c(t) = (x^i(t),u^{\alpha}(t))$, the $\rho$-admissibility of $c$ is expressed by 
the relation $\dot{x}^i(t)= \gamma^i_{\alpha}(x(t))u^{\alpha}(t)$ for all $t \in [0,1]$
(see Section 2). Next, take any point $e_0 = (x^i_0,y^A_0) \in E_{\nu(n_0)}$ and 
consider the following system of linear first-order ordinary differential equations with 
time-dependent coefficients:
\[
\dot{y}^A = \Gamma^A_{\alpha B}(x(t))u^{\alpha}(t)y^B.
\]
It follows from the theory of linear differential equations that this system admits a
unique solution $y^A(t)$ with $y^A(0) = y^A_0$ and which, moreover, is defined for all
$t \in [0,1]$. The curve $c^h(t) = (x^i(t),y^A(t))$ then clearly satisfies all the 
requirements of the proposition.

The proof for the more general case, with $\hbox{Im}(c)$ not necessarily contained 
in a single bundle chart, follows by taking a partition $0 = t_0 < t_1 < ... <t_n = 1$
of $[0,1]$ in such a way that the previous construction can be applied to the restriction
of $c$ to each subinterval $[t_i,t_{i+1}]$, and then glueing the results together.
\QED 

We will call $c^h$ the {\it $h$-lift} of the admissible curve $c$, with initial
point $e_0$. In coordinates it follows from the above that an $h$-lift of a 
$\rho$-admissible curve $c(t) = (x^i(t),u^{\alpha}(t))$ in $N$ is a curve 
$c^h(t) = (x^i(t),y^A(t))$ in $E$ for which
\begin{equation}\label{hlift}
\dot{x}^i(t) = \gamma^i_{\alpha}(x(t))u^{\alpha}(t), \quad 
\dot{y}^A(t)= \Gamma^A_{\alpha B}(x(t))u^{\alpha}(t)y^B(t)\;.
\end{equation} 
It can be immediately inferred from these relations that, in general, $c^h$ is
not fully determined by the projection $\tilde{c}(t) = (x^i(t))$ of $c$ alone. 
More precisely, different $\rho$-admissible curves projecting onto the same base 
curve in $M$ may have different $h$-lifts in $E$ with the same initial point.

Proposition \ref{liftcurve} allows us to associate a notion of parallel transport
to a linear $\rho$-connection. Indeed, consider a smooth admissible curve
$c: [0,1] \rightarrow N$ with projection $\tilde{c} = \nu \circ c$ on $M$ and put
$\tilde{c}(0) = m_0, \tilde{c}(1) = m_1$. One can then define a map
\[
\tau_c: E_{m_0} \longrightarrow E_{m_1},\; e_0 \longmapsto c^h(1),
\]
where $c^h$ is the $h$-lift of $c$ with initial point $c^h(0) = e_0$. From the 
construction of the $h$-lift it easily follows that this map is indeed well-defined 
and, moreover, determines a linear isomorphism between the fibres $E_{m_0}$ and 
$E_{m_1}$. We will call $\tau_c$ {\it the operator of parallel transport (or parallel
displacement) along the $\rho$-admissible curve $c$}. It is important to emphasise
again that, in general, parallel transport can not be unambiguously associated to
a base curve in $M$.

The construction of $\tau_c$ can obviously be extended to the case where $c$ is 
a piecewise smooth admissible curve. In order to introduce a suitable concept of 
holonomy in the framework of linear $\rho$-connections, it turns out that the class 
of admissible curves in $N$ should be further extended to curves admitting 
(a finite number of) discontinuities in the form of certain `jumps' in the fibres 
of $N$, such that the corresponding base curve is piecewise smooth. 
A detailed discussion of this matter will be the topic of a separate paper. For a 
treatment of holonomy in the special case where $(N,\rho)$ defines a 
Lie-algebroid structure on $M$: see, for instance, the recent papers by Fernandes 
\cite{Fern1,Fern2}.

In the next subsection, we describe the construction of an operator which for 
linear $\rho$-connections can be seen as the analogue of the covariant derivative 
operator in standard connection theory. 

\subsection{The associated derivative operator}
Consider a linear $\rho$-connection $h$ on the vector bundle 
$\pi$, with associated connection map $K$ (\ref{K1}). Take $s \in \Gamma(\nu)$ and
$\psi \in \Gamma(\pi)$. For any $m \in \hbox{Dom}(s) \cap \hbox{Dom}(\psi)$ one
readily verifies that $(s(m),{\psi}_{\ast}(\rho(s(m))))$ determines an element of 
the bundle $\rho^{\ast}TE$. We then define $\nabla_s\psi \in \Gamma(\pi)$ by
\begin{equation}\label{nabla1}
\nabla_s\psi(m) = K(s(m),{\psi}_{\ast}(\rho(s(m))))\,.   
\end{equation}
Let $U \subset \hbox{Dom}(s) \cap \hbox{Dom}(\psi)$ be a trivialising coordinate 
neigbourhood for both $\nu$ and $\pi$, with coordinates $x^i$ on $U$ and 
corresponding local bundle coordinates $(x^i,u^{\alpha})$ and $(x^i,y^A)$ on $N$ 
and $E$, respectively. Putting $s(x) = (x^i,s^{\alpha}(x)), \; \psi(x) = 
(x^i,\psi^A(x))$, we then find, using (\ref{K2}):
\begin{equation}\label{nabla2}
\nabla_s\psi(x) = 
\left(x^i, \frac{\partial {\psi}^A}{\partial x^j}(x)\gamma^j_{\alpha}(x)s^{\alpha}(x) - 
\Gamma^A_{{\alpha}B}(x)s^{\alpha}(x)\psi^B(x)\right)\,.    
\end{equation}
In terms of the vector field $X = \rho \circ s \in {\frak X}(M)$, we can still rewrite
the components of $\nabla_s\psi$ as
\[
(\nabla_s\psi)^A(x) = 
\frac{\partial {\psi}^A}{\partial x^j}(x)X^j(x) - 
\Gamma^A_{{\alpha}B}(x)s^{\alpha}(x)\psi^B(x)\,.
\] 
The following theorem gives a full characterisation of the operator $\nabla$, whereby
it is tacitly assumed that its action is restricted to those pairs 
$(s,\psi) \in \Gamma(\nu) \times \Gamma(\pi)$ for which $\hbox{Dom}(s)$ and
$\hbox{Dom}(\psi)$ have nonempty intersection.

\begin{theorem}\label{deriv} The operator 
$\nabla: \Gamma(\nu) \times \Gamma(\pi) \rightarrow \Gamma(\pi)$, 
defined by (\ref{nabla1}), satisfies the following properties:\\
(i) $\nabla$ is $\R$-bilinear;\\
(ii) for all $(s,\psi) \in \Gamma(\nu) \times \Gamma(\pi)$ and $f \in C^{\infty}(M)$
we have:
\[
\nabla_{fs}\psi = f\nabla_s\psi \qquad \hbox{and} \qquad \nabla_s(f\psi) = 
f\nabla_s\psi + (\rho \circ s)(f)\psi\,.
\]
Moreover, $\nabla$ is uniquely determined by the given linear $\rho$-connection $h$.
\end{theorem}
{\bf Proof.} The proofs of the properties (i) and (ii) follow by straightforward
computation. The fact that $\nabla$ is uniquely determined by $h$ can be easily
deduced from (\ref{nabla1}) and the definition of the connection map $K$. Indeed,
different $\rho$-connections necessarily induce different maps $V$ (see (\ref{V})) and, 
hence, different connection maps $K$ (see (\ref{K1})).
\QED

We will call the operator $\nabla$ the {\it derivative operator associated
to the linear $\rho$-connection $h$}. In case $N=TM$ and $\rho$ is the identity map
on $TM$, we recover the classical notion of covariant derivative operator of a linear 
connection on a vector bundle over $M$. In his treatment of Lie algebroid connections 
on a vector bundle, where $N=A$ is a Lie algebroid over $M$ with anchor 
map $\rho$, Fernandes refers to the $\nabla$-operator as the $A$-derivative: 
see \cite{Fern2}.
 
From the fact that $\nabla_s\psi$ is $C^{\infty}(M)$-linear in $s$, it follows that
for a given $\psi$, $(\nabla_s\psi)(m)$ only depends on the value of $s$ in $m$, and 
not on the behaviour of $s$ in a neighbourhood of $m$. This allows us to define 
for each $n \in N$, with $m = \nu(n)$, an operator
\begin{equation}\label{nabla3}
\nabla_n: \Gamma_m(\pi) \longrightarrow E_m, \psi \longmapsto \nabla_n\psi := 
\nabla_s\psi(m),
\end{equation}
where $s$ may be any (local) section of $\nu$ for which $s(m) = n$. Alternatively,
we could have defined defined the operator $\nabla_n$ directly according to the 
prescription $\nabla_n\psi = K(n,\psi_{\ast}(\rho(n))$. The properties
of $\nabla_n$ immediately follow from Theorem \ref{deriv}, i.e. $\nabla_n$ is
$\R$-linear and for any $f \in C^{\infty}(M)$ and $\psi \in \Gamma_m(\pi)$,
we have that 
\[
\nabla_n(f\psi) = f(m)\nabla_n\psi + \rho(n)(f)\psi(m).
\]
Next, let $c: I \rightarrow N$ be an admissible curve in $N$, with corresponding 
base curve $\tilde{c}= \nu \circ c$. Consider a map $\tilde{\psi}: I \rightarrow E$,
i.e.\ a curve in $E$, satifying $\pi \circ \tilde{\psi} = \tilde{c}$. It is now 
readily seen that, for each $t \in I$, 
$(c(t),\dot{\tilde{\psi}}(t)) \in \rho^{\ast}TE$ and we may then define
\[
\nabla_c\tilde{\psi}(t):= K(c(t),\dot{\tilde{\psi}}(t))\;,
\]
which we will call {\it the derivative of $\tilde{\psi}$ along the admissible curve $c$}. 
In coordinates, putting $c(t) = (\tilde{c}^i(t),c^{\alpha}(t))$ and $\tilde{\psi}(t) = 
(\tilde{c}^i(t), \tilde{\psi}^A(t))$, we obtain
\[
(\nabla_c\tilde{\psi}(t))^A = \frac{d\tilde{\psi}^A}{dt}(t) - 
\Gamma^A_{\alpha B}(\tilde{c}(t))c^{\alpha}(t)\tilde{\psi}^B(t)\;.
\]
Assume one can find a (local) section $\psi \in \Gamma(\pi)$ such that 
$\psi(\tilde{c}(t)) = \tilde{\psi}(t)$ for all $t \in I$. This will be the case, 
for instance, if the base curve $\tilde{c}$ is an injective immersion. A straightforward 
computation then shows that 
\begin{eqnarray*}
(\nabla_{c(t)}\psi)^A &=& \frac{\partial \psi^A}{\partial x^j}\dot{\tilde{c}}^j(t)
- \Gamma^A_{\alpha B}c^{\alpha}(t)\psi^B(\tilde{c}(t))\\
&=& \frac{d\tilde{\psi}^A}{dt}(t) - 
\Gamma^A_{\alpha B}(\tilde{c}(t))c^{\alpha}(t)\tilde{\psi}^B(t)\;,
\end{eqnarray*}
where, for the second equality, we have used the fact that $\psi^A(\tilde{c}(t))
\equiv \tilde{\psi}^A(t)$.
We may therefore conclude that the derivative of $\tilde{\psi}$ along $c$ verifies
\[
\nabla_c\tilde{\psi}(t) = \nabla_{c(t)}\psi,
\]
for any $\psi \in \Gamma(\pi)$ such that $\psi(\tilde{c}(t)) \equiv \tilde{\psi}(t)$,
if such a section $\psi$ exists. 

\begin{remark} {\rm A special situation occurs when $\rho$ has a nontrivial kernel and  
the image of an admissible curve $c$ is contained in it. In particular, we then know 
that $c(t)$ necessarily belongs to a fixed fibre of $\nu$ (cf. Section 2) and 
the base curve $\tilde{c}$ reduces to a point in $M$, say 
$\tilde{c}(t) = \nu(c(t)) = m_0$ for all $t$. We then consider a map 
$\tilde{\psi}: I \rightarrow E_{m_0}$. In coordinates, 
with $m_0 = (x^i_0), \tilde{\psi}(t) =(x^i_0,y^A(t))$, we then find that 
\[
\nabla_c\tilde{\psi}(t) = \left(x^i_0, \dot{y}^A(t)- \Gamma^A_{\alpha B}(x_0)c^{\alpha}(t)y^B(t)\right)\in E_{m_0}\;.
\]  
In particular, if we associate to each point $e_0 = (x^i_0,y^A_0) \in E_{m_0}$  
the constant map $\tilde{\psi}(t) \equiv (x^i_0,y^A_0)$, we 
obtain a time-dependent linear map on the fibre $E_{m_0}$, namely 
$e_0 \mapsto \nabla_ce_0(t) = (x^i_0, - \Gamma^A_{\alpha B}(x_0)c^{\alpha}(t)y^B_0)$}.    
\end{remark}

Next, it is easy to see how the action of the derivative 
operator of a linear $\rho$-connection on a vector bundle $\pi:E \rightarrow M$, 
can be extended to sections of the dual vector bundle $\pi^{\ast}: E^{\ast} 
\rightarrow M$. If, by convention, for $s \in \Gamma(\nu)$ and $f \in C^{\infty}(M)$ 
we put $\nabla_sf = (\rho \circ s)(f)$, we can immediately define an action of the 
operator $\nabla_s$ on $\Gamma(\pi^{\ast})$ as follows: for any ${\frak f} \in 
\Gamma(\pi^{\ast})$, $\nabla_s{\frak f} \in \Gamma(\pi^{\ast})$ is uniquely determined 
by
\[
\langle \psi,\nabla_s{\frak f}\rangle = \nabla_s \langle \psi,{\frak f}\rangle 
- \langle \nabla_s\psi,{\frak f}\rangle\,,
\]
for all $\psi \in \Gamma(\pi)$, where $\langle\;,\;\rangle$ denotes the canonical
pairing between sections of $\pi$ and sections of $\pi^{\ast}$. Herewith, it is then
standard to further extend the action of $\nabla_s$ to sections of any 
tensor bundle constructed out of $E$ and $E^{\ast}$.

In what precedes we have shown that a linear $\rho$-connection on a vector bundle
$\pi: E \rightarrow M$ gives rise to an operator $\nabla$ verifying the conditions of
Theorem \ref{deriv}. We now demonstrate that the converse also holds.

\begin{theorem} Any operator $\nabla: \Gamma(\nu) \times \Gamma(\pi) \rightarrow
\Gamma(\pi)$, verifying the properties (i) and (ii) of Theorem \ref{deriv}, is the
derivative operator of a unique linear $\rho$-connection on $\pi$ .
\end{theorem}
{\bf Proof.} Take $n \in N$, with $\nu(n)= m$, and $\psi \in \Gamma_m(\pi)$. From the
above discussion it follows that the given operator $\nabla$ induces an operator
$\nabla_n$ on $\Gamma_m(\pi)$ such that $\nabla_n\psi \in E_m$. Putting $\psi(m) = e$, 
and denoting by $\iota_e: E_m \rightarrow V_eE$ the canonical isomorphism 
between the vector spaces $E_m$ and $V_eE$, we may consider 
the vector $\psi_{\ast}(\rho(n)) - \iota_e(\nabla_n\psi) \in T_eE$. It is now 
straightforward to check that the mapping $\Gamma_m(\pi) \rightarrow TE,
\psi \mapsto \psi_{\ast}(\rho(n)) - \iota_e(\nabla_n\psi)$ is $C^{\infty}(M)$-linear 
in $\psi$ and, hence, only depends on the value of $\psi$ in $m$. From this we deduce
that there exists a well-defined smooth mapping $h: \pi^{\ast}N \rightarrow TE$, 
given by 
\[
h(e,n) = \psi_{\ast}(\rho(n)) - \iota_e(\nabla_n\psi)\,,
\]
for any $\psi \in \Gamma(\pi)$ with $\psi(\nu(n)) = e$. Clearly, $\pi_{\ast}(h(e,n))
= \rho(n)$, which already shows that $\rho \circ \tilde{\pi}_2 = \pi_{\ast} \circ h$.
The linearity of $h_e = h(e,.): N_{\pi(e)} \rightarrow T_eE$ is obvious.    
With $\{\phi_t\}$ denoting the flow of the dilation vector field on $E$, and observing
that for any $\psi \in \Gamma(\pi)$ we also have $\phi_t \circ \psi \in \Gamma(\pi)$
for each $t \in \R$, it is not difficult to verify that 
\begin{eqnarray*}
(\phi_t)_{\ast}\left(\psi_{\ast}(\rho(n)) - \iota_e(\nabla_n\psi)\right) &=&
(\phi_t \circ \psi)_{\ast}(\rho(n)) - \iota_{\phi_t(e)}(\nabla_n(\phi_t \circ \psi))\\
&=&h(\phi_t(e),n)\,,
\end{eqnarray*}
proving that $h$ is indeed a linear $\rho$-connection.

It now remains to be shown that the given $\nabla$ is the derivative operator of the 
constructed $\rho$-connection $h$. Let $K$ denote the connection map associated to $h$
(cf. Section 3). Using (\ref{K1}) and (\ref{V}), together with the above definition
of $h$, we find for any $n \in N$ and $\psi \in \Gamma_{\nu(n)}(\pi)$, putting 
$\psi(\nu(n)) = e$:
\begin{eqnarray*}
K(n,\psi_{\ast}(\rho(n)))&=& p_2(\psi_{\ast}(\rho(n)) - h(e,n))\\
&=&p_2(\iota_e(\nabla_n\psi))\\
&=&\nabla_n\psi\,.
\end{eqnarray*}
Since, in view of (\ref{nabla1}), the left-hand side precisely determines the 
derivative operator associated to $h$, this completes the proof of the theorem.
\QED  

Suppose $\nabla$ and $\bar{\nabla}$ are the derivative operators corresponding
to two (different) linear $\rho$-connections on the vector bundle $\pi$. It follows
from Theorem \ref{deriv} that the difference $\nabla - \bar{\nabla}$ is a 
$C^{\infty}(M)$-bilinear mapping $S: \Gamma(\nu) \times \Gamma(\pi) 
\rightarrow \Gamma(\pi)$, which locally reads
\[
(S(s,\psi))^A = (\Gamma^A_{\alpha B} - \bar{\Gamma}^A_{\alpha B})s^{\alpha}\psi^B\;.
\]
Conversely, given a derivative operator $\nabla$ and an arbitrary 
$C^{\infty}(M)$-bilinear mapping $S: \Gamma(\nu) \times \Gamma(\pi) 
\rightarrow \Gamma(\pi)$, the operator $\nabla + S$, mapping any pair
of sections $(s,\psi)$ onto $\nabla_s\psi + S(s,\psi)$, also defines a 
derivative operator verifying the assumptions of Theorem \ref{deriv} and, hence, 
determines a linear $\rho$-connection on $\pi$. We also note that $S$ uniquely 
determines a smooth section ${\cal S}$
of the tensor product bundle $N^{\ast}\otimes E^{\ast} \otimes E \rightarrow M$, 
where $N^{\ast} \rightarrow M$ and $E^{\ast} \rightarrow M$ are the dual 
bundles of $N \rightarrow M$ and $E \rightarrow M$, respectively. The relation
between $S$ and $\cal S$ is given by
\[
{\cal S}(m)(n,e,e^{\ast}) = \langle S(s,\psi)(m),e^{\ast}\rangle\;,
\]
for all $m \in M, n \in N_m, e \in E_m, e^{\ast}\in E^{\ast}_m$,
and where $s$ and $\psi$ are any sections of $\nu$ and $\pi$, respectively,
such that $s(m) = n$ and $\psi(m) = e$. 

\section{Curvature and torsion}
Clearly, in the case of arbitrary vector bundles 
$\nu: N \rightarrow M$ and $\pi: E \rightarrow M$ there is no way, in general, of 
assigning a notion of torsion or curvature to a linear $\rho$-connection. 
However, let us assume in what follows that the space 
of sections $\Gamma(\nu)$ is equipped with an algebra structure (over $\R$), 
with product denoted by $*$, such that the mapping $\Gamma(\nu) \times \Gamma(\nu)
\rightarrow \Gamma(\nu), (s_1,s_2) \mapsto s_1*s_2$
is $\R$-bilinear and skew-symmetric and, in addition, verifies a Leibniz-type
rule
\begin{equation}\label{Leibniz}
s_1*(fs_2) = f(s_1*s_2) + \rho(s_1)(f)s_2\;,
\end{equation}
for all $s_1,s_2 \in \Gamma(\nu)$ and $f \in C^{\infty}(M)$.
Note that we do not require $\rho$ to induce an algebra morphism between 
$(\Gamma(\nu),*)$ and $({\frak X}(M),[\;,\;])$.

Whenever the space of sections of the vector bundle $\nu: N \rightarrow M$ 
is equipped with a bilinear operation $*$ verifying the above assumptions, we will
follow \cite{Grab1} in saying that $N$ admits the structure of a 
{\it pre-Lie algebroid\/}. When dropping the skew-symmetry assumption of the
product $*$, we obtain a so-called pseudo-Lie algebroid (cf. \cite{Grab1}).
For a treatment of the differential calculus on pseudo- and pre-Lie algebroids, we
refer to \cite{Grab2}, where both structures are simply called ``algebroids" and
``skew algebroids", respectively. The algebraic counterpart of pre-Lie algebroids,
namely differential pre-Lie algebras, have also been studied in \cite{Kos}.
  
In analogy with the Poisson structure that exists on the dual bundle of any Lie
algebroid, one can show that on the dual bundle $\mu: N^{\ast} \rightarrow M$
of any pre-Lie algebroid $\nu: N \rightarrow M$ there exists a distinguished
bivector field $\Lambda$ which, in particular, induces an `almost-Poisson' bracket 
$\{\;,\;\}$ on $C^{\infty}(N^{\ast})$, verifying all properties of a Poisson bracket 
except for the Jacobi identity. One can show that the Schouten-Nijenhuis bracket
of the bivector field $\Lambda$ with itself vanishes (and, hence, $\Lambda$ becomes
a Poisson tensor) iff the algebra $(\Gamma(\nu),*)$ is a Lie algebra, i.e.\
the Jacobi identity holds for the $*$-product. In that case one can also prove 
that $\rho$ induces a Lie algebra homomorphism between $(\Gamma(\nu),*)$ and 
$({\frak X}(M),[\;,\;])$, and $N$ then becomes a Lie algebroid over $M$.
   
\subsection{Curvature}
Assume $N$ admits a pre-Lie algebroid structure, with product $*$ on $\Gamma(\nu)$,
and consider a linear $\rho$-connection on a vector bundle $\pi: E \rightarrow M$, 
with associated derivative operator $\nabla$. We may now define a mapping 
$R: \Gamma(\nu) \times \Gamma(\nu) \times \Gamma(\pi) \rightarrow \Gamma(\pi)$
given by
\begin{equation}\label{curvature}
R(s_1,s_2;\psi) := \nabla_{s_1}\nabla_{s_2}\psi - \nabla_{s_2}\nabla_{s_1}\psi
- \nabla_{s_1*s_2}\psi\;.
\end{equation}
It easily follows that $R$ is $C^{\infty}(M)$-linear and skew-symmetric in 
$s_1$ and $s_2$, but fails to be $C^{\infty}(M)$-linear in $\psi$.
Indeed, a straightforward computation shows that for arbitrary $s_1,s_2 \in 
\Gamma(\nu), \psi \in \Gamma(\pi)$ and $f \in C^{\infty}(M)$,
\[
R(s_1,s_2;f\psi) = fR(s_1,s_2;\psi) + \left(\rho(s_1)\circ \rho(s_2) - 
\rho(s_2)\circ \rho(s_1) - \rho(s_1*s_2)\right)(f)\psi\;.
\]
From this expression it is seen that $R$ will be fully tensorial iff 
$\rho$ induces an algebra homomorphism from $(\Gamma(\nu),*)$ to 
$({\frak X}(M),[\;,\;])$,
i.e.
\begin{equation}\label{lie}
\rho(s_1*s_2) = [\rho(s_1),\rho(s_2)].
\end{equation}
In particular, this implies that for all $s_1,s_2,s_3 \in \Gamma(\nu)$
we have
\[
s_1*(s_2*s_3) + s_2*(s_3*s_1) + s_3*(s_1*s_2) \in \Gamma(\nu_{|\ker (\rho)})\;, 
\]
i.e. the `Jacobiator' of the $*$-product should take values in the kernel of
the vector bundle morphism $\rho$. (The denomination `Jacobiator' is sometimes
used in the literature to indicate, in an algebra with a skew-symmetric
product, the cyclic sum that vanishes in case the Jacobi identity holds). 
If (\ref{lie}) holds, we will call the mapping $R$, defined by (\ref{curvature}), 
the {\it curvature\/} of the given $\rho$-connection. 

\begin{remark}{\rm Another important consequence of (\ref{lie}) is that the 
generalised distribution ${\cal D} (= \hbox{Im}(\rho))$ on $M$ is involutive. 
Note, however, that since $\rho$ need not be of constant rank, involutivity does 
not necessarily imply integrability of ${\cal D}$. (For integability conditions 
of a generalised distribution, see e.g.\ \cite{Suss,Vais1}.)} 
\end{remark} 

Consider a local coordinate neighbourhood $U$ in $M$, with coordinates $x^i\;
(i = 1,\ldots,n)$, which is also a trivialising neighbourhood for 
both vector bundles $\nu$ and $\pi$. Let $\sigma_{\alpha}\; (\alpha = 1, \ldots, k)$, 
respectively $p_A\; (A = 1, \ldots, \ell)$, represent a local basis of sections 
of $\nu$, respectively $\pi$, defined on $U$. We then have
\[
\sigma_{\alpha}*\sigma_{\beta} = c^{\lambda}_{\alpha \beta}\sigma_{\lambda}\;,
\]
for some functions $c^{\lambda}_{\alpha \beta} \in C^{\infty}(U)$.
Putting $\rho(\sigma_{\alpha}) = \gamma_{\alpha}^{i}(\partial/\partial x^i)$,
the condition (\ref{lie}) yields the following relation
\[
c^{\lambda}_{\alpha \beta}\gamma^i_{\lambda}
= \gamma^j_{\alpha}\frac{\partial \gamma^i_{\beta}}{\partial x^j}
- \gamma^j_{\beta}\frac{\partial \gamma^i_{\alpha}}{\partial x^j}\;,
\]
for all $\alpha,\beta,i$. Given a linear $\rho$-connection on $\pi$, let
\[
\nabla_{\sigma_{\alpha}}p_A = \Gamma^B_{\alpha A}p_B\;.
\]
Denoting the components of the curvature $R$ with respect to the chosen
local bases of sections by $R^B_{\alpha \beta A}$, i.e.\  
$R(\sigma_{\alpha},\sigma_{\beta};p_A) = R^B_{\alpha \beta A}p_B$, a
straightforward computation reveals that
\begin{equation}\label{compon}
R^B_{\alpha \beta A} = \gamma^i_{\alpha}\frac{\partial \Gamma^B_{\beta A}}
{\partial x^i} - \gamma^i_{\beta}\frac{\partial \Gamma^B_{\alpha A}}
{\partial x^i} + \Gamma^B_{\alpha C}\Gamma^C_{\beta A} -
\Gamma^B_{\beta C}\Gamma^C_{\alpha A} - 
c^{\lambda}_{\alpha \beta}\Gamma^B_{\lambda A}\;.
\end{equation}

Always under the assumption that (\ref{lie}) is satisfied, we will establish 
a link between the curvature of a linear $\rho$-connection $h$ on 
$\pi: E \rightarrow M$ and the (lack of) involutivity of the (generalised)
distribution ${\cal Q} = \hbox{Im}$$(h)$. Recalling that for any $s \in \Gamma(\nu)$, 
$s^h \in {\frak X}(E)$ denotes its $h$-lift (cf. Section 2), we have the 
following useful property.

\begin{lemma} For any $s_1,s_2 \in \Gamma(\nu)$
\[
[s_1^h,s_2^h](e) - (s_1*s_2)^h(e) \in V_eE \quad \hbox{for all}\; e \in E.
\]
\end{lemma}
{\bf Proof.} From the fact that for each $s \in \Gamma(\nu)$, $s^h$ and 
$\rho \circ s$ are $\pi$-related vector fields (cf.\ Proposition \ref{sh} (iii)),
it follows that $[s_1^h,s_2^h]$ and $[\rho(s_1),\rho(s_2)]$ are also
$\pi$-related. Taking into account (\ref{lie}) we then easily find that
\begin{eqnarray*}
\pi_{\ast}\left([s_1^h,s_2^h] - (s_1*s_2)^h\right)&=&\pi_{\ast}[s_1^h,s_2^h]
- \rho(s_1 * s_2)\circ \pi\\
&=&\left([\rho(s_1),\rho(s_2)] - \rho(s_1*s_2)\right)\circ \pi \\
&=&0\;,
\end{eqnarray*}
from which the result follows.
\QED

We now come to the following important result which tells us that the curvature 
$R$ can indeed be seen as a measure for the `non-involutivity' of the (generalised) 
distribution $\cal Q$ determined by a linear $\rho$-connection. 
(Recall that $\iota_e$ denotes the canonical identification between 
$E_{\pi(e)}$ and $V_eE$).

\begin{theorem} For any $s_1, s_2 \in \Gamma(\nu)$ we have
\begin{equation}\label{integrability}
\iota_e\left([s_1^h,s_2^h](e) - (s_1*s_2)^h(e)\right) = R(s_1,s_2;\psi)(m)\;,
\end{equation}
for each $e \in E$ for which the left-hand side is defined, and
where $m = \pi(e)$ and $\psi$ is any section of $\pi$ such that $\psi(m) = e$.
\end{theorem}
{\bf Proof.} First of all, note that the left-hand side of (\ref{integrability})
makes sense in view of the previous lemma, and that the `tensorial character' of
$R$ implies that the right-hand side does not depend on the choice of the 
section $\psi$ for which $\psi(m) =e$. Secondly, using the properties of
the $h$-lift of sections it is not difficult to check that
$[s_1^h,s_2^h] - (s_1*s_2)^h$ is $C^{\infty}(M)$-linear
in both $s_1$ and $s_2$. Since we already know that the same is true for 
$R(s_1,s_2;\psi)$, it suffices to verify (\ref{integrability}) on a local basis 
of sections $(\sigma_{\alpha})_{\alpha = 1,\ldots ,n}$ of $\Gamma(\nu)$, defined 
on a suitable coordinate neighbourhood $U$ of $m = \pi(e)$, for some chosen point
$e \in E$. There is no loss of generality by assuming that $U$ is also a 
trivialising neighbourhood for $\pi$, and denote the corresponding bundle 
coordinates on $E$ by $(x^i,y^A)$. 
In particular, let the coordinates of the point $e$ be given by $(x^i_0,y^A_0)$.

Using the notations introduced above, we find after a rather tedious, but
straightforward computation, that
\[
[\sigma_{\alpha}^h,\sigma_{\beta}^h](e) - (\sigma_{\alpha}*\sigma_{\beta})^h(e)
=\left(\gamma^i_{\alpha}\frac{\partial \Gamma^A_{\beta B}}{\partial x^i}
- \gamma^i_{\beta}\frac{\partial \Gamma^A_{\alpha B}}{\partial x^i}
+ \Gamma^A_{\beta C}\Gamma^C_{\alpha B} - 
\Gamma^A_{\beta C}\Gamma^C_{\alpha B}
- c^{\lambda}_{\alpha \beta}\Gamma^A_{\lambda B}\right)_{x_0}y^B_0
\frac{\partial}{\partial y^A}_{|_e}\;.
\] 
The result now easily follows when comparing the right-hand side with the expression 
(\ref{compon}) for the local components of $R$, and bearing in mind that $\iota_e$ maps 
each $(x^i_0,y^A_0,0,w^A) \in V_eE$ onto $(x^i_0,w^A) \in E_m$.
\QED

{\bf Example.} If $(N,\nu)$ is a Lie algebroid over $M$ with anchor map $\rho$,
we recover the notion of curvature defined, for instance, in \cite{Fern2}.

\subsection{Torsion}
Assume again $\nu:N \rightarrow M$ is a pre-Lie algebroid, i.e.\ that 
$\Gamma(\nu)$ admits an algebra structure, with a skew-symmetric product $*$ satisfying (\ref{Leibniz}). We do not require, however, that $\rho$
is an algebra homomorphism. Consider now a linear $\rho$-connection on $\nu$, 
with associated derivative operator $\nabla$ (i.e.\ we take $E = N$ and $\pi = \nu$). 
We can then define a mapping $T: \Gamma(\nu) \times \Gamma(\nu) \rightarrow 
\Gamma(\nu)$ given by
\begin{equation}\label{torsion}
T(s_1,s_2) = \nabla_{s_1}s_2 - \nabla_{s_2}s_1 - s_1*s_2\;.
\end{equation}
It is not difficult to check that $T$, which may be called the {\it torsion\/}
of the given $\rho$-connection, is a $C^{\infty}(M)$-bilinear and 
skew-symmetric mapping. Let $(\sigma_{\alpha})_{\alpha = 1, \ldots, k}$ represent
a local basis of sections of $\nu$ such that 
\[
\nabla_{\sigma_{\alpha}}\sigma_{\beta} = \Gamma^{\lambda}_{\alpha \beta}\sigma_{\lambda}
\quad \hbox{and} \quad \sigma_{\alpha}*\sigma_{\beta} =
c^{\lambda}_{\alpha \beta}\sigma_{\lambda}\;. 
\]
 
It then readily follows that
\[
T(\sigma_{\alpha},\sigma_{\beta}) = (\Gamma^{\lambda}_{\alpha \beta} -  
\Gamma^{\lambda}_{\beta \alpha} - c^{\lambda}_{\alpha \beta})\sigma_{\lambda}\;.
\]

{\bf Example.} Let $A$ be a type $(1,1)$-tensor field on $M$ and consider a linear 
pseudo-connection on $\tau_M$ with fundamental tensor field $A$ (cf.\ Section 4 (iv)).
Here we have $N = TM$, $\nu = \tau_M$ and for the product $*$ we may take the bracket 
$[\;,\;]_A$ on ${\frak X}(M)$, defined by (\ref{Nijenhuis}). This bracket satisfies 
(\ref{Leibniz}), but in general will not be a Lie bracket (since $A$ need not be a
Nijenhuis tensor). The notion of torsion, defined by (\ref{torsion}), corresponds to 
the one encountered in treatments of pseudo-connections (see e.g.\ \cite{Etayo,Wong}). 

\section{Principal $\rho$-connections}

As before, let $\nu: N \rightarrow M$ be a vector bundle and 
$\rho: N \rightarrow TM$ a vector bundle map, such that $\tau_M \circ \rho = \nu$.
Let $\pi: P \rightarrow M$ be a principal $G$-bundle, with a free  
(say, right) group action $\Phi: P \times G \rightarrow P$, such that $P/G \cong M$. 
The Lie algebra of $G$ will be denoted by $\frak g$. For a standard treatment of 
the theory of principal bundles, we refer to \cite{KoNo}. Using again the notations 
$\Phi(e,g) = \Phi_g(e) = eg$, recall from Section 2 that a principal $\rho$-connection 
on a principal $G$-bundle $\pi: P \rightarrow M$ is a $\rho$-connection $h$ on 
$\pi$ satisfying the additional condition $(\Phi_g)_{\ast}(h(e,n))=h(eg,n)$, 
for all $g \in G$ and $(e,n) \in \pi^{\ast}N$, i.e.\ $h$ is equivariant with respect 
to the induced actions of $G$ on $\pi^{\ast}N$ and $TP$. 

In this section we will briefly describe some aspects of the theory
of principal $\rho$-connections. Much more on the subject can be found, for instance, 
in \cite{Fern2} for the case where $(N,\nu)$ is a Lie algebroid over $M$ with anchor 
map $\rho$. In fact, all properties described in that paper which do not effectively
rely on the Lie algebra structure of $\Gamma(\nu)$, also hold in the more general 
setting we are considering here. 

First of all, given a principal $\rho$-connection $h$ on $\pi$, the 
$G$-equivariance of $h: \pi^{\ast}N \rightarrow TP$ implies that it induces
a bundle mapping from $\pi^{\ast}N/G$ into $TP/G$. Taking into account
that $\pi^{\ast}N/G \cong N$, and putting $\hat{\pi}: TP \rightarrow TP/G$ the 
natural projection onto the space of $G$-orbits, we obtain a well-defined mapping 
\[
\bar{\omega}: N \longrightarrow TP/G\;, n \longmapsto \hat{\pi}(h(e,n))\;, 
\]
for any $e$ such that $\pi(e) = \nu(n)$.  

At this point it is important to recall that $TP/G$ admits a Lie algebroid 
structure over $M$ (for details, see Appendix A of \cite{Mackenzie}). In particular, 
we have a vector bundle structure 
\[
\hat{\tau}: TP/G \rightarrow M
\] 
such that the following diagram commutes

\begin{picture}(15,4.6)(-3,4.5)
\thicklines
\put(2.4,5.4){$P$}
\put(3.2,5.6){\vector(1,0){3}}
\put(6.5,5.4){$M$}
\put(2.7,7.6){\vector(0,-1){1.7}}
\put(6.5,8){$TP/G$}
\put(3.2,8.2){\vector(1,0){3}}
\put(2.2,8){$TP$}
\put(6.8,7.6){\vector(0,-1){1.7}}
\put(4.7,5.2){$\pi$}
\put(4.7,8.4){$\hat{\pi}$}
\put(7,6.7){$\hat{\tau}$}
\put(2.1,6.9){$\tau_P$}
\end{picture}     

and the $C^{\infty}(M)$-module of sections $\Gamma(\hat{\tau})$ is equipped with a 
Lie bracket which we shall denote here by $[\;,\;]^{\wedge}$. We also recall that,
given a local trivialising neigbourhood $U \subset M$ of the principal bundle $\pi$, 
we have the identification $T_{U}P/G \cong TU \times {\frak g}$. The anchor map
${\frak p}: TP/G \rightarrow TM$ of the Lie algebroid structure on $TP/G$ 
precisely corresponds to the projection onto the first factor in this local splitting.

Following Fernandes \cite{Fern2}, let $\Gamma^p(N^{\ast},TP/G) := 
\Gamma(\nu^{{\ast}(p)})\otimes \Gamma(\hat{\tau})$ denote the $C^{\infty}(M)$-module 
of $TP/G$-valued sections of the exterior bundle $\nu^{{\ast}(p)}:\bigwedge^pN^{\ast} \rightarrow M$, where $\nu^{\ast}: N^{\ast} \rightarrow M$ is the dual bundle of $N$. 
Clearly, the mapping $\bar{\omega}$, defined above, is a vector bundle
mapping over the identity on $M$ and, hence, we can associate to it
a unique element $\omega \in \Gamma^1(N^{\ast},TP/G)$ according to $\omega(s):= 
\bar{\omega}\circ s$, for arbitrary $s \in \Gamma(\nu)$. We will call $\omega$ the 
{\it connection $1$-section\/} of the given principal $\rho$-connection $h$. From the
previous definitions it can be easily deduced that
\begin{equation}\label{compatible}
{\frak p} \circ \omega = \rho\;.
\end{equation}
In addition, we have the following interesting property, the proof of which 
is also an immediate consequence of the definitions of the various objects involved.

\begin{proposition} Given a principal $\rho$-connection $h$ on $\pi$, with associated
connection $1$-section $\omega$, then for any $s \in \Gamma(\nu)$
\[
\omega(s) \circ \pi = \hat{\pi} \circ s^h\;.
\]
\end{proposition}

Assume now that $\Gamma(\nu)$ is equipped with an algebra structure, with skew-symmetric
bilinear product $*$ satisfying (\ref{Leibniz}), such that $N$ becomes a pre-Lie 
algebroid. We may then put, for arbitrary $s_1, s_2 \in \Gamma(\nu)$,
\[
\Omega(s_1,s_2):= [\omega(s_1),\omega(s_2)]^{\wedge} - \omega(s_1 * s_2)\;.
\]
Taking into account (\ref{compatible}) one immediately verifies that $\Omega$ is $C^{\infty}(M)$-bilinear, and since it is obviously skew-symmetric, we have that 
$\Omega$ is an element of $\Gamma^2(N^{\ast},TP/G)$, called the {\it curvature
$2$-section} of the principal $\rho$-connection.  

Let $U \subset M$ be a local trivialising neighbourhood of the principal 
bundle $\pi$. Then, given a principal $\rho$-connection 
on $\pi$ with associated connection $1$-section $\omega$, the isomorphism 
$T_UP/G \cong TU \times {\frak g}$ allows one to write, for any local section
$s \in \Gamma(\nu)$ defined on $U$, 
\[
\omega(s) = (\rho(s),\omega_U(s))\;.
\]
This uniquely determines an element $\omega_U \in \Gamma(N^{\ast}_{|U},{\frak g})$,
called a {\it local connection $1$-section}.
When considering an open covering $\{U_j\}$ of $M$ by trivialising neighbourhoods 
of $\pi$, one can associate in this way a local connection $1$-section 
$\omega_j \equiv \omega_{U_j}$ to each $U_j$. Moreover, for any two overlapping
neighbourhoods $U_j$ and $U_k$, with corresponding transition function $\psi_{jk}:
U_j \cap U_k \rightarrow G$ (cf. \cite{KoNo}, Vol.1), one can show that the relation
\begin{equation}\label{gauge}
\omega_k = \hbox{Ad}(\psi_{jk}^{-1})\omega_j + \psi_{jk}^{-1}d_{\rho}(\psi_{jk})
\end{equation}
holds on $U_j \cap U_k$, with $\psi^{-1}_{jk}(m):= (\psi_{jk}(m))^{-1}$. 
Here, $\hbox{Ad}$ denotes the adjoint representation
of $G$ on $\frak g$, whereas $d_{\rho}$ is an operator which associates
to any smooth mapping $f: M \rightarrow G$, the mapping
\[
d_{\rho}f: N \longrightarrow TG,\; n \longmapsto f_{\ast}(\rho(n))\;.
\]
(The definition of $d_{\rho}$ can be extended to maps from $M$ into any smooth 
manifold: see e.g. \cite{Fern2}). Conversely, given an open covering $\{U_j\}$ of 
$M$ by trivialising neighbourhoods and any family of $\frak g$-valued $1$-sections
$\omega_j$ (each defined on the corresponding $U_j$) for which (\ref{gauge}) holds,
one can demonstrate that there exists a unique principal $\rho$-connection $h$ on
$\pi$ for which the $\omega_j$'s are local connection $1$-sections.
For a proof of these results, and for more details on local
connection $1$-sections as well as on the notion of local curvature 
$2$-section, associated to a principal $\rho$-connection, we refer to 
\cite{Fern1,Fern2}. 

Finally, it is not difficult to see that also in the present framework, a principal 
$\rho$-connection on a principal $G$-bundle $\pi: P \rightarrow M$, induces a 
$\rho$-connection on any fibre bundle associated to $P$ (cf.\ Vol.1 of \cite{KoNo}
for the construction in the standard case, and \cite{Fern2} for the Lie 
algebroid case).

\section{Final remarks} 
In this paper we have described a general framework for connections on fibre
bundles $\pi: E \rightarrow M$, defined over a vector bundle map 
$\rho: N \rightarrow TM$, with $\nu: N \rightarrow M$ a given vector bundle. Our 
main source of inspiration was provided by some recent work of R.L. Fernandes,
who treated the case where $(N,\nu)$ is a Lie algebroid with anchor map $\rho$ 
\cite{Fern2}. By dropping the requirement that $N$ should be equipped with a
Lie algebroid structure, we have extended the picture to the case where
the distribution $\rho(N)$ on $M$  need not be integrable. In that
respect, one of us (BL) has found some interesting applications of the 
theory in sub-Riemannian geometry \cite{Bavo1}, 
as well as in a new connection theoretic approach to nonholonomic mechanics 
\cite{Bavo2}. Further work along these lines, also in the field of geometric control
theory, is in progress.  
Also from a purely geometrical point of view, there are several aspects of
the theory which still need further investigation. In particular, we intend 
to study in more detail the notion of parallel transport and the concept of holonomy 
in the general setting described above. Moreover, in analogy with Fernandes' treatment 
of contravariant connections in Poisson geometry, it may also be of interest to
investigate the role of connections over a bundle map induced by some other 
geometric structures on a manifold (cf. Section 4).
 
While finalising the present paper, we came across a preprint of a recent paper 
by M. Popescu and P. Popescu \cite{Pop1}.
From this paper we learned that some of the ideas developed above are probably 
closely related to work done by these authors in the past decade. In particular,
the idea of a generalised $\rho$-connection on a vector bundle seems to be contained 
in a paper from 1992, entitled ``On the geometry of relative tangent spaces" \cite{Pop2}. 
A relative tangent space, called an `anchored vector bundle' in \cite{Pop1}, 
precisely refers to a structure consisting of a vector bundle 
$\nu: N \rightarrow M$ and an `anchor map' $\rho: N \rightarrow TM$. 

{\bf Acknowledgements}
This work has been partially supported by a research grant from the ``Bijzonder
Onderzoekfonds" of Ghent University (Belgium).

\end{document}